\documentclass{amsart}
\usepackage{amsfonts,amssymb,amscd,amsmath,enumerate,verbatim,calc}
\usepackage[all]{xy}
\SelectTips{eu}{}

\newcommand{\fn}{\mathfrak{n}}
\newcommand{\fm}{\mathfrak{m} }
\newcommand{\fq}{\mathfrak{q} }

\newcommand{\BN}{\mathbb N}
\newcommand{\BZ}{\mathbb Z}

\newcommand{\lra}{\longrightarrow}
\newcommand{\xra}{\xrightarrow}

\newcommand{\pd}{\operatorname{projdim}}
\newcommand{\id}{\operatorname{injdim}}

\newcommand{\rank}{\operatorname{rank}}

\newcommand{\Ann}{\operatorname{Ann}}

\newcommand{\Ker}{\operatorname{Ker}}

\newcommand{\preg}{\operatorname{polyreg}}
\newcommand{\post}{\operatorname{post}}

\newcommand{\avind}{\boldsymbol{A}}
\newcommand{\lind}{\boldsymbol{L}}

\newcommand{\depth}{\operatorname{depth}}
\newcommand{\edim}{\operatorname{embdim}}

\newcommand{\syz}[3]{\Omega^{#2}_{#1}(#3)}

\newcommand{\Tor}{\operatorname{Tor}}
\newcommand{\tor}[4]{\operatorname{Tor}_{#1}^{#2}(#3,#4)}
\newcommand{\Hom}{\operatorname{Hom}}
\newcommand{\hh}{\operatorname{H}}
\newcommand{\Ext}{\operatorname{Ext}}
\newcommand{\ext}[4]{\operatorname{Ext}^{#1}_{#2}(#3,#4)}

\newcommand{\grass}[2]{\operatorname{gr}_{#1}(#2)}
\newcommand{\ges}{\geqslant}

\newcommand{\hilb}[2]{\operatorname{Hilb}_{#1}(#2)}
\newcommand{\hiltop}[2]{h_{#2}(#1)}
\newcommand{\hsp}[2]{{H\!S}_{#1}(#2)}
\newcommand{\hspoly}[4]{{\tau}_{#1}^{#2}(#3;#4)}
\newcommand{\length}{{\ell}}

\theoremstyle{plain}

\theoremstyle{plain}

\newtheorem{theorem}{Theorem}[section]

\newtheorem{proposition}[theorem]{Proposition}
\newtheorem{lemma}[theorem]{Lemma}
\newtheorem*{Lemma}{Lemma}
\newtheorem{corollary}[theorem]{Corollary}

\newtheorem{itheorem}{Theorem}

\newtheorem{icorollary}[itheorem]{Corollary}
\newtheorem{iproposition}[itheorem]{Proposition}

\theoremstyle{definition}

\newtheorem{chunk}[theorem]{}

\theoremstyle{remark}
\newtheorem*{Question}{Question}
\newtheorem{remark}[theorem]{Remark}

\numberwithin{equation}{theorem}

\begin{document}
\title[Hilbert-Samuel functions]
{Hilbert-Samuel functions of modules over Cohen-Macaulay rings}

\author{Srikanth Iyengar} 
\address{Department of Mathematics, University of Nebraska,  Lincoln, NE 68588, U.S.A.}  
\email{iyengar@math.unl.edu}

\author{Tony J. Puthenpurakal} 
\address{Department of Mathematics, IIT Bombay, Powai, Mumbai 400 076, India} 
\email{tputhen@math.iitb.ac.in}

\thanks{S.I. was partly supported by NSF grant DMS 0442242}

\subjclass{Primary 13D40; Secondary 13D02, 13D07}

\date{\today}
 
 \begin{abstract}
   For a finitely generated, non-free module $M$ over a CM local ring $(R,\fm,k)$, it is
   proved that for $n\gg 0$ the length of $\tor 1RM{R/\fm^{n+1}}$ is given by a polynomial
   of degree $\dim R-1$.  The vanishing of $\tor iRM{N/\fm^{n+1}N}$ is studied, with a
   view towards answering the question: if there exists a finitely generated $R$-module
   $N$ with $\dim N\ge 1$ such that the projective dimension or the injective dimension of
   $N/\fm^{n+1}N$ is finite, then is $R$-regular? Upper bounds are provided for $n$ beyond
   which the question has an affirmative answer.
\end{abstract}

\maketitle

\section*{Introduction}
Let $(R, \fm,k)$ denote a local ring with maximal ideal $\fm$ and residue field $k$; in this
article local rings are assumed to be Noetherian.  Let $M$ be a finite, that is to say,
finitely generated $R$-module.

Recall that the \emph{Hilbert-Samuel function} of $M$ is the function on the non-negative
integers that maps $n$ to $\length_R(M/\fm^{n+1}M)$, where $\length_R(-)$ denotes length.
The Hilbert-Samuel function is an important invariant of $M$, and has for long been a
topic of active research.  Note that $M/\fm^{n+1}M=\tor 0RM{R/\fm^{n+1}}$. This suggests
the following line of enquiry: Fix an integer $i\ge0$ and consider the function 
\[
n\mapsto \length_R \, \tor iRM{R/\fm^{n+1}} \quad \text{for each} \quad n\in \BN
\]
This may be thought of as the \emph{$i$th Hilbert-Samuel function of $M$}.  For each $i$, this
function is given by a polynomial for $n\gg 0$, which we denote $\hspoly iRMz$. This
result is classical for $i=0$, and the polynomial $\hspoly 0RMz$ is the
\emph{Hilbert-Samuel} polynomial of $M$. For general $i$, it is due to Kodiyalam
\cite{VKov}; see also Theodorescu \cite{ET}.

Recall that $\deg \hspoly 0RMz=\dim M$.  Our main result, contained in the theorem below,
is a lower bound on the degree of $\hspoly iRMz$ for $i\geq 1$. The upper bound is already
in \cite[(2)]{VKov}, and is valid even when $\fm$ is substituted by any $\fm$-primary
ideal. A convention: the degree of a polynomial $t(z)$ is $-1$ if and only if $t(z)=0$.

\begin{itheorem}
\label{igrowth}
Let $R$ be a local ring with $\depth R\geq 1$, and $M$ a non-zero
finite $R$-module. For $i\in\BN$, if $\pd_RM\geq i$, then $\tor iRM{R/\fm^{n+1}}\ne 0$ for
each integer $n\ge0$, and
\[
\dim R - 1 \geq \deg\hspoly iRMz \geq \depth R - 1
\]
\end{itheorem}

Specializing to the case where $R$ is Cohen-Macaulay yields the following result; it
subsumes \cite[(18)]{Pu1}, which assumes in addition that $M$ is maximal Cohen-Macaulay.

\begin{icorollary}
\label{cmgrowth}
Let $R$ be a Cohen-Macaulay local ring with $\dim R\geq 1$, and $M$ a non-zero
finite $R$-module. For $i\in\BN$, if $\pd_RM\geq i$, then
\[
\deg\hspoly iRMz =\dim R - 1
\]
\end{icorollary}

This contrasts drastically with the general situation:

\begin{iproposition}
\label{example:noncm}
Given non-negative integers $p,q$ with $p\leq q-1$, there exists a local ring $R$ with
$\depth R=p+1$ and $\dim R = q$, and a finite maximal Cohen-Macaulay $R$-module $M$ with
$\deg\hspoly 1RMz=p$.
\end{iproposition}

Therefore, in Corollary \ref{cmgrowth} it is crucial that $R$ is Cohen-Macaulay.  However,
we have been unable to find a module $M$ for which the inequalities in Theorem
\ref{igrowth} are strict.  On the other hand, towards the end of Section \ref{growth} we
provide an example that shows that it is not possible to replace $\fm$ with an arbitrary
$\fm$-primary ideal.

The preceding results are proved in Section \ref{growth}. Our proof of the lower bound in
Theorem \ref{igrowth}, that $\deg\hspoly iRMz\geq \depth R -1$, proceeds via induction on
$\depth R$.  The non-vanishing of Tor is an immediate consequence of the following
elementary observation, which appears as Lemma \eqref{minor lemma} in the text:

\smallskip

\emph{If $\tor 1RM{R/\fm^{n+1}}=0$ for some integer $n\ge0$, then $\fm^n\syz 1RM=0$; hence
  either $M$ is free or $\depth R=0$.}

\smallskip

In what follows $\syz dRM$ denotes the $d$th syzygy in the minimal free resolution of $M$.
In Theorem \eqref{Intheorem}, the remark above is extended to statements concerning the
vanishing of $\tor iRM{N/\fm^{n+1}N}$, where $N$ is another finite $R$-module. These are
akin to a result of Levin and Vasconcelos~\cite[(1.1)]{LV}, stated as Theorem
\eqref{theorem:lv}. Any one of these, applied with $M=k$ and $N=R$, may be used to deduce
that if either $\pd_R(R/\fm^{n+1})$ or $\id_R(R/\fm^{n+1})$ is finite, then $\fm^{n+1}=0$
or $R$ is regular.

The discussion in the preceding paragraph suggests the question: If there is a finite
$R$-module $N$ with $\dim_RN\geq 1$ and an integer $n\ge0$ such that either
$\pd_R(N/\fm^{n+1}N)$ or $\id_R(N/\fm^{n+1}N)$ is finite, then is $R$ regular?  Note that
one cannot expect such a strong conclusion if one drops the restriction on the dimension
of $N$.  Indeed, over Cohen-Macaulay rings (with dualizing modules) it is easy to construct
modules of finite length and finite projective (injective) dimension.

Section 3 focuses on this question. The main result is stated below; $\rho_R(N)$ is
defined in \eqref{rho vs polyreg}, while $\avind_R(N)$ is defined in
\eqref{avind:definition}.  These numbers are bounded above by the polynomial regularity of
$N$, and so finite; see \eqref{rho vs polyreg} and Lemma \eqref{avind}.

\begin{itheorem}
\label{iext}
Let $(R,\fm,k)$ be a local ring and let $N$ be a finite $R$-module with $\depth_RN\ge 1$
and $\id_{R}(N/\fm^{n+1}N)$ finite for some non-negative integer $n$.  If $n \geq
\rho_R(N)$, then $R$ is a hypersurface; if $n\geq \avind_R(N)$ as well, then $R$ is
regular.
\end{itheorem}

This is the content of Theorem \eqref{nth2}; an analogue for projective dimension is
given in Theorem \eqref{nth}.  Our proof of the result above uses superficial sequences,
and has a different flavour from that of Theorem \eqref{nth}, which is easily deduced from
available literature. Superficial sequences and other techniques traditionally used in the study
of Hilbert functions play an important part in our arguments in Sections 2 and 3.  The
relevant definitions and results are recalled in Section 1.

\section{Preliminaries}
Let $(R,\fm,k)$ be a local ring, with maximal ideal $\fm$ and residue field $k$. Let $M$
be a finite $R$-module. We write $\grass \fm M$ for $\bigoplus_{n\ges 0}\fm^nM/\fm^{n+1}M$
and view it as a graded module over the associated graded ring $\grass \fm R$.

\begin{chunk}
\label{hilbert}
Recall that the \emph{Hilbert function} of $M$ assigns to each integer $n$ the value
$\rank_k(\fm^n M/\fm^{n+1}M)$. Its generating series is the \emph{Hilbert series} of $M$;
we denote it $\hilb Mz$. This series can be expressed as a rational
function of the form
\[
\hilb Mz  = \frac{\hiltop zM}{(1-z)^{\dim M}} \quad \text{with} \quad  \hiltop zM \in \BZ[z]
\]
The \emph{multiplicity} of $M$ is the integer $e(M) = \hiltop 1M$; see, \cite[(4.6)]{BH} for
details.
\end{chunk}

\begin{chunk}
\label{superficial}
  Let $x$ be an element in $\fm$. When $x$ is non-zero, let $n$ be the largest integer
  such that $x\in \fm^n$, and let $x^*$ denote the image of $x$ in $\fm^n/\fm^{n+1}$.
  Set $0^* = 0$.
  
  The element $x \in \fm$ is \emph{superficial} for $M$ if there exists a positive
  integer $c$ with
  $$
  (\fm^{n+1}M \colon_M x)\cap\fm^cM = \fm^{n}M \ \text{ for all }\quad n \geq c
  $$
\end{chunk}

The following properties of superficial elements are often invoked in this work.

\begin{chunk}
\label{superficial=regular}
Let $M_1,\dots,M_s$ be a finite $R$-modules.
\begin{enumerate}[\quad\rm(a)]
\item 
If $k$ is infinite, then there exists an element $x$ that is superficial for each
  $M_j$.
\item If an element $x$ is superficial on $M_j$ and $\depth_RM_j\geq 1$, then $x$ is a
  non-zero divisor on $M_j$ and $(\fm^{n+1}M_j\colon_{M_j} x) = \fm^{n}M_j$ for $n \gg 0$.
\end{enumerate}

Indeed, it suffices to verify the desired properties for $M=M_1\oplus\cdots \oplus M_s$.
Now (a) is trivial when $\dim M=0$.  When $\dim M\geq 1$, an element $x$ is superficial on
$M$ if and only if $x\notin \fm^2$ and $x^*$ is not in any associated prime ideal of
$\grass \fm M$, except perhaps $\grass \fm R_{\ges 1}$. Thus, superficial elements exist
if the residue field $k$ is infinite.  As to (b), an argument of Sally \cite[p.\ 7]{Sa}
for $M=R$ extends to the general case.
\end{chunk}

Lead by property (b) above, for each element $x$ superficial on $M$ we set
\[
\rho_R(x,M) =  \inf\{r \mid (\fm^{n+1}M\colon_M x) = \fm^{n}M\quad \text{for each $n\geq r$}\}
\]

We provide an upper bound for this number; it involves the following invariant.

\begin{chunk}
\label{polyreg}
Let $\mathcal M$ denote the irrelevant maximal ideal of $\grass \fm R$, and $\hh_{\mathcal
  M}^j(\grass \fm M)$ the $j$th local cohomology of $\grass\fm M$ with respect to
$\mathcal M$. Set
\[
\preg_R(M) = \sup\{p \mid \hh_{\mathcal M}^j(\grass \fm M)_{p-j}\ne 0
       \quad \text{for some $j$}\}
\]
This number is  called the \emph{polynomial regularity} of $M$ over $R$; see
\c{S}ega \cite[(1.6)]{Sega}.
\end{chunk}

\begin{lemma}
\label{rho:lemma}
If $\depth_RM\geq 1$, then $\rho_R(x,M)\leq \preg_R(M)+1$.
\end{lemma}

\begin{remark}
\label{rho vs polyreg}
Thanks to the lemma, the number
\[
\rho_R(M)=\sup\{\rho_R(x,M)\mid \text{$x$ a element superficial on $M$}\}
\]
satisfies the inequality $\rho_R(M)\leq\preg_R(M)+1$, and is in particular finite.  This
result is subsumed in \cite{Pu3}, where it is proved that $\rho_R(x,M)$ is independent of
$x$ and the bound is improved to $\rho_R(x,M)\leq \preg_R(M)$.  So we provide only a

\begin{proof}[Sketch of the proof of Lemma \eqref{rho:lemma}]
One has an equality of formal Laurent series
\[
\sum_{n\geq 0}\length_R(M/\fm^{n+1}M)\, z^n= \frac{\hilb Mz}{(1-z)}\,;
\]
this is immediate from \eqref{hilbert}; the formal power series on the left is the
Hilbert-Samuel series of $M$. In particular, there exists an integer $c$ and a polynomial
$\hsp Mz$ such that $\length_R(M/\fm^{n+1}M)=\hsp Mn$ for $n\geq c$.  The least such
number $c$ is called the postulation number of $M$; we denote it $\post_R(M)$.  Set
$N=M/xM$. Arguing along the lines of the proof of Elias' result \cite[(1.1)]{Elsc}, one
obtains
\[
\rho_R(x,M)\leq \max\{\post_R(M), \post_R(N)\} + 1
\]
Now $\post_R(M)\leq \preg_R(M)$ and $\post_R(N)\leq \preg_R(N)$; this is immediate from
\cite[(4.4.3)]{BH}, and the definition of $\preg_R(-)$ recalled above.  Moreover, $x$ is a
superficial non-zero divisor on $M$, by \eqref{superficial=regular}, so $\preg_R(N)\leq
\preg_R(M)$. Combining the inequalities above yields the desired result.
\end{proof}
\end{remark}

A standard trick allows one to assume that superficial elements exist:

\begin{chunk}
\label{infresfield}
Let $R[X]$ be a polynomial ring over $R$ on a finite set $X$ of variables.  Set $R' =
R[X]_{\fm'}$, where $\fm'= \fm R[X]$, and $M'= M \otimes_R R'$.  The ring $R'$ is again
local, with maximal ideal $\fm R'$, which we again denote $\fm'$, and residue field $k' =
k(X)$, the field of rational functions over $k$. The following claims are verified easily.
\begin{enumerate}[\quad\rm (a)]
\item
 $\length_{R} M = \length_{R'}M'$.
\item
$\hspoly iRMz = \hspoly i{R'}{M'}z$ for each non-negative integer $i$.
\item 
$\edim R = \edim {R'}$, $\dim M = \dim M'$ and $\depth M = \depth M'$.
\item
$\id_R(M/\fm^{n}M)$ is finite if and only if $\id_{R'} M'/(\fm')^{n}M'$ is finite.
\end{enumerate}
\end{chunk}

\section{Growth}
\label{growth}
In this section we prove Theorem \ref{igrowth} from the introduction.  Parts of the
arguments are abstracted out in the following lemmas.  The one below is extended in
Theorem \eqref{Intheorem}; the crux of the argument is simple and well illustrated in this
special case.

 \begin{lemma}
 \label{minor lemma}
 Let $(R,\fm,k)$ be a local ring and $M$ a finite $R$-module.  If there is an integer
 $n\ge0$ with $\tor 1RM{R/\fm^{n+1}}=0$, then $\fm^n\syz 1RM=0$; hence either $M$ is free
 or $\depth R=0$.
 \end{lemma}

 \begin{proof}
   Set $L=\syz 1RM$ and let $0\to L \to F\to M\to 0$ be the exact sequence defining $L$.
   Applying $-\otimes_R{R/\fm^{n+1}}$ to it yields the exact sequence
 \[
 0\to L/\fm^{n+1}L \to F/\fm^{n+1}F \to M/\fm^{n+1}M\to 0
 \]
 Thus $\fm^{n+1}L = L\cap \fm^{n+1}F$. By construction, $L\subseteq \fm F$ so
 $\fm^{n}L\subseteq L\cap \fm^{n+1}F$, hence $\fm^{n}L\subseteq \fm^{n+1}L$, so Nakayama's
 lemma yields: $\fm^n L =0$.  It remains to note that $L$ is a submodule of the free module
 $F$.
 \end{proof}
 
 Next we recall the following result, contained in \cite[(2)]{VKov}.

\begin{chunk}
\label{form}
Let $(R,\fm,k)$ be a local ring, $M$ a finite $R$-module, and $i$ a positive integer.
The function on the non-negative integers defined by
\[
n\mapsto \length_R\, \tor iRM{R/\fm^{n+1}R}
\]
is given by a polynomial for $n \gg 0$, and of degree at most $\dim R - 1$.
\end{chunk}


\begin{remark}
\label{shift}
We let $\hspoly iRMz$ denote the polynomial that arises in the \eqref{form}.  Since $\tor
iRM-\cong \tor 1R{\syz{i-1}RM}-$ for $i\geq 1$, it follows that
\[
\hspoly iRMz \cong \hspoly 1R{\syz {i-1}RM}z \quad\text{for each} \quad i\geq 1
\]
This equality often allows one to obtain results on $\hspoly iRMz$ from corresponding
statements concerning the case $i=1$.
\end{remark}
 
 Remark \eqref{shift} focuses our attention on $\hspoly 1RMz$; the next result describes
 its basic properties. Recall our convention that the degree of the polynomial $0$ is
 $-1$.

\begin{lemma}
\label{hspoly:properties}
Let $(R,\fm,k)$ be a local ring and $M$ a finite $R$-module.
\begin{enumerate}[\quad\rm(1)]
\item
For each finite $R$-module $N$, one has 
\[
\hspoly 1R{M\oplus N}z = \hspoly 1RMz + \hspoly 1RNz\,,
\]
in particular, $\deg\hspoly 1R{M\oplus N}z=\sup\{\deg\hspoly 1RMz,\deg\hspoly 1RNz\}$.
\item 
If $\depth R\ge 1$ and $M$ is not free, then $\deg \hspoly 1RMz\ge 0$.
\item 
If $0\to L\to M\to N\to 0$ is an exact sequence of $R$-modules with $\length_R L$
finite, then $\deg\hspoly 1RMz\ge \deg\hspoly 1RNz$.
\item
If $\length_RM$ is finite, then $\deg\hspoly 1RMz=\dim R -1$.
\end{enumerate}
\end{lemma}

\begin{proof}
  The first assertion is immediate from the isomorphisms
\[
\tor 1R{M\oplus N}{R/\fm^{n+1}} \cong \tor 1RM{R/\fm^{n+1}} \oplus \tor 1RN{R/\fm^{n+1}}
\]
while (2) is contained in Lemma \eqref{minor lemma}.

\smallskip

Proof of (3): For each non-negative integer $n$, applying $-\otimes_R{R/\fm^{n+1}}$ to the
given exact sequence yields an exact sequence of $R$-modules
\[
\cdots \to \tor 1RM{R/\fm^{n+1}}\to \tor 1RN{R/\fm^{n+1}}\to L/\fm^{n+1}L \to M/\fm^{n+1}M
\to \cdots
\]
Since $\length_RL$ is finite, the map $L/\fm^{n+1}L\to M/\fm^{n+1}M$ is injective for $n\gg 0$.
For such $n$, computing lengths in the sequence above yields
\[
\length_R\, \tor1RM{R/\fm^{n+1}} \ge \length_R\, \tor1RN{R/\fm^{n+1}}
\]
This implies the desired result.

\smallskip

Proof of (4): It suffices to prove that $\deg\hspoly 1RMz\ge\dim R-1$; this is due to
\eqref{form}.  In view of part (3) above, a standard induction on the length
reduces the claim to the case $M=k$. Note that for each non-negative integer $n$ one has
an isomorphism of $R$-modules
\[
\tor 1Rk{R/\fm^{n+1}}\cong (\fm\cap \fm^{n+1})/\fm^{n+2}=\fm^{n+1}/\fm^{m+2}
\]
Thus, the degree of $\hspoly 1RMz$ equals the degree of the Hilbert polynomial of $R$,
that is to say, to $\dim R-1$. The settles (3).
\end{proof}

The proof of Theorem \ref{igrowth} is an induction on $\depth R$, and uses
the following

\begin{lemma}
\label{igrowth:lemma}
Let $(R,\fm,k)$ be a local ring and $M$ a finite non-free $R$-module.  Let $x\in \fm$ be a
superficial non-zero divisor on $R$, $M$, and $\syz 1RM$.  Set $S=R/xR$ and $N=M/xM$.
For $n\gg 0$, one has
\begin{align*}
&\hspoly 1RMn = \hspoly 1RM{n-1} + \hspoly 1SN{n-1} \\
&\deg \hspoly 1RMz = \deg \hspoly 1SNz  + 1
\end{align*}
\end{lemma}

\begin{proof}
  Let $\fn$ be the maximal ideal of $S$. By \eqref{superficial=regular}, since $x$ is
  superficial on $R$, for each $n \gg 0$ one has an exact sequence of
  $R$-modules
\[
0 \to R/{\fm^n} \xra{\kappa} {R}/{\fm^{n+1}} \to {S}/{\fn^{n+1}} \to 0
\]
where $\kappa$ maps the residue class of $a$ in $R/\fm^n$ to the residue class of $xa$ in
$R/\fm^{n+1}$, for each $a\in R$. Applying $M\otimes_R-$ to it yields an exact sequence of
$R$-modules
\begin{align*}
  \tor 1RM{R/\fm^{n}} 
    &\xra{\tor 1RM{\kappa}} \tor 1RM{R/\fm^{n+1}} 
              \to \tor 1RM{S/\fn^{n+1}} \to  \\
\to {M}/{\fm^{n}M} &\xra{M\otimes_R\kappa} {M}/{\fm^{n+1}M} \to  {N}/{\fn^{n+1}N} \to 0
\end{align*}
The maps $M\otimes_R\kappa$ and $\tor 1RM{\kappa}$ are injective for $n\gg0$. Indeed, the
former is injective by \eqref{superficial=regular}, because $x$ is superficial on $M$. As
to the latter, set $L=\syz 1RM$ and consider the commutative diagram
\[
\xymatrixrowsep{2pc} 
\xymatrixcolsep{2.5pc} 
\xymatrix{ 
0\ar@{->}[r] &\tor 1RM{R/\fm^{n}}\ar@{->}[r]\ar@{->}[d]^{\tor 1RM{\kappa}}
                 & L/\fm^n L \ar@{->}[d]^{L\otimes_R\kappa} \\
0\ar@{->}[r] &\tor 1RM{R/\fm^{n+1}}\ar@{->}[r] & L/\fm^{n+1}L
}
\]
obtained from the exact sequence $0\to L \to F\to M\to 0$ defining $L$ by applying
$-\otimes_R{R/\fm^n}$ and $-\otimes_R{R/\fm^{n+1}}$. By choice, $x$ is superficial on $L$,
so for $n\gg 0$ the map $L\otimes_R\kappa$ is injective, again by
\eqref{superficial=regular}. Thus, $\tor 1RM{\kappa}$ is also injective.

Returning to the long exact sequence above with the injectivity on hand, one finds
for each $n\gg 0$ the following exact sequence of $R$-modules 
\[
0\to \Tor^{R}_{1}(M,R/\fm^{n}) \to \tor 1RM{R/\fm^{n+1}}
\to \Tor^{R}_{1}(M , S/\fn^{n+1})\to 0
\]
Now $\Tor^R_1(M,S/\fn^{n+1})\cong \Tor^{S}_{1}(N,S/\fn^{n+1})$, as $x$ is both $R$-regular
and $M$-regular.  This isomorphism and the exact sequence above yield, for $n\gg0$, an
equality
\[
 \length_R\, \tor 1RM{R/\fm^{n+1}} =
            \length_R\, \Tor^{R}_{1}(M,R/\fm^{n}) + \length_R\, \Tor^{S}_{1}(N,S/\fn^{n+1}) 
\]
Thus, $\hspoly 1RMn = \hspoly 1RM{n-1} + \hspoly 1SN{n}$. The hypothesis implies $\syz
1RM\ne 0$ and $\depth R\ge 1$, so $\tor 1RM{R/\fm^{n+1}}$ is non-zero for each $n\ge 0$,
by Lemma \eqref{minor lemma}. Consequently, the displayed equality yield $\deg \hspoly
1RMz = \deg \hspoly 1SNz + 1$, even when $\tor iSN{S/\fn^{n+1}}=0$ for $n\gg 0$, that is to say,
when $\deg\hspoly 1SNz=-1$.
\end{proof}

\begin{chunk}{\emph{Proof of Theorem \ref{igrowth}}.}
\label{igrowth:proof}
In view of Remark \eqref{shift}, it suffices to verify the theorem for $i=1$; the
hypothesis on $M$ is then that it is not free.  The first part of the theorem: that $\tor
1RM{R/\fm^{n+1}}\ne 0$ for each $n\ge 0$, is immediate from Lemma \eqref{minor lemma}.
Given \eqref{form}, it remains to prove that $\deg\hspoly 1RMz\geq \depth R -1 $.

To establish this we induce on $\depth R$; the base case $\depth R=1$ is easily settled:
use Lemma (\ref{hspoly:properties}.2).  Assume that the desired inequality has been
verified for local rings of depth $d$, for some $d\ge 1$.  Let $R$ be a local ring of
depth $d+1$ and $M$ a finite, non-free $R$-module.

When $\dim_RM=0$, the desired result is provided by Lemma (\ref{hspoly:properties}.4).

Suppose $\dim_RM\ge 1$.  By Lemma (\ref{hspoly:properties}.1) we may assume $M$ is
indecomposable. Set $L = \{m\in M\mid \fm^s\cdot m=0\quad \text{for some $s\geq 1$}\}$.
Since the length of $L$ is finite it suffices to verify the inequality for $M/L$, by Lemma
(\ref{hspoly:properties}.3).  Now $\depth_R(M/L)\ge 1$ and $M/L$ is not free; the last
property holds because $M$ is indecomposable and not free. Thus, substituting $M/L$ for
$M$ one may assume $\depth_RM\ge 1$.

By \eqref{infresfield}, one can extend the residue field of $R$ and ensure that it is
infinite, so by \eqref{superficial=regular} there is an element $x$ that is a superficial
non-zero divisor on $R$, $M$, and $\syz 1RM$.  Hence the local ring $S=R/xR$ has depth
$d$, and the $S$-module $N=M/xM$ is not free, since $\pd_SN=\pd_RM$.  Lemma
\eqref{igrowth:lemma} now provides the equality below
\[
\deg \hspoly 1RMz = \deg\hspoly 1SNz + 1 \geq \depth S - 1 + 1 = \depth R - 1
\]
while the induction hypothesis yields the inequality. \qed
\end{chunk}

Next we prepare for the proof of Proposition \ref{example:noncm}.

\begin{chunk}
\label{trivial extensions}
Let $(S,\fn,k)$ be a local ring and let $L$ be a finitely generated $S$-module. Set
$R=S\ltimes L$; thus $R=S\oplus L$, with multiplication given by $(s,l)\cdot (s',l') =
(ss', sl'+s'l)$.  Evidently, $R$ is a local ring with maximal ideal $\fm = \fn\oplus L$.
We view $S$ as an $R$-module via the canonical surjection $R\to S$. In particular
\begin{equation*}
\dim_RS = \dim S \qquad\text{and}\qquad \depth_RS = \depth S 
\end{equation*}

Note that $\dim R=\dim S$ and $\depth R = \min\{\depth S,\depth_S L\}$, so $R$ is
Cohen-Macaulay if and only if $S$ is Cohen-Macaulay and the $S$-module $L$ is maximal
Cohen-Macaulay.  The information of interest to us is contained in the following

\begin{Lemma}
  For each non-negative integer $n$ one has
\[
\length_R\, \tor 1RS{R/\fm^{n+1}} = \rank_k (\fn^{n}L/\fn^{n+1}L)
\]
In particular, $\deg\hspoly 1RSz=\dim_SL - 1$.
\end{Lemma}
\begin{proof}
Note that the kernel of the canonical surjection $R\to S$ is $L$, viewed as an
ideal of $R$. This implies the first isomorphism below:
\[
\tor 1RS{R/\fm^{n+1}}\cong (L \cap \fm^{n+1})/\fm^{n+1} L \cong \fn^{n}L/\fn^{n+1}L 
\]
The second equality holds because $\fm = \fn \oplus L$ and $L^2=0$. Now compute lengths.
\end{proof}
\end{chunk}


\begin{proof}[Proof of Proposition \ref{example:noncm}]
  Let $k$ be a field, $S=k[[x_1,\dots,x_{q}]]$ the ring of formal power series in
  variables $x_1,\dots,x_{q}$, and let $L=S/(x_{p+2},\dots,x_q)$. The local ring $S$ is
  Cohen-Macaulay with $\dim S = q$ and $L$ is a finite Cohen-Macaulay $S$-module with
  $\dim_SL=p+1$. Set $R=S\ltimes L$, and let $M=S$, viewed as an $R$-module via the
  canonical surjection $R\to S$.  By \eqref{trivial extensions}, one has $\depth R = p+1$,
  $\dim R = q$, $\deg\hspoly 1RMz=p$, and $M$ is a maximal Cohen-Macaulay $R$-module.
\end{proof}

Another interpretation of this proposition is as follows: for any integer $s\ge 1$ there
exists a local ring $R$ and a maximal Cohen-Macaulay $R$-module $M$ such that $\dim R -
\deg\hspoly 1RMz=s$. Thus, the inequality in \cite[(4.a)]{ET} is strict, in general.

  It is also interesting to study the growth of the function defined by
\[
n\mapsto \length_R\,\tor iRM{R/I^{n+1}}
\]
for an arbitrary $\fm$-primary ideal $I$. In this case as well the function is given by a
polynomial for $n\gg 0$, and its degree is at most $\dim R - 1$, by \cite[(2)]{VKov}.
However, upper bound may be strict even if $R$ is Cohen-Macaulay: when the residue field
of such an $R$ is infinite and $M$ is a non-free maximal Cohen-Macaulay module, there are
$\fm$-primary ideals $I$ such that $\tor 1RM{R/I^{n+1}}=0$ for $n\gg 0$; see
\cite[(20)]{Pu1}.

\section{Finite homological dimensions}
\label{homological dimensions}
The results in this section are motivated by the

\begin{Question}
\label{quest}
Let $R$ be a local ring and $N$ a finite $R$-module such that, for some integer $n\geq0$,
the $R$-module $N/\fm^{n+1}N$ has either finite projective dimension or finite injective
dimension.  If $\dim_RN\geq 1$, then is $R$ regular?
\end{Question}

These questions are suggested by Theorem \eqref{Intheorem} below and its corollaries.  As
noted in the introduction, some restriction on $\dim_RN$ or on $\depth_RN$ is crucial.
We begin with the result below; with $N=R$ and $n=1$ it captures Lemma \eqref{minor lemma}.


\begin{theorem}
\label{Intheorem}
Let $(R,\fm,k)$ be a local ring, and let $M$ and $N$ be finite $R$-modules.  If there
exist non-negative integers $i$ and $n$ such that
\[
\tor iRMN = 0 = \tor iRM{N/\fm^{n+1}N}
\]
then $\fm^n(\syz iRM \otimes_R N)=0$ and either $\depth_R N=0$ or $\pd_RM\leq i-1$.
\end{theorem}

With stronger hypotheses, we are able to arrive at a better conclusion.

\begin{theorem}
\label{Intheorem2}
Let $(R,\fm,k)$ be a local ring, and let $M$ and $N$ be finite $R$-modules. 
If there exist non-negative integers $i$ and $n$ such that
\[
\tor iRM{N/\fm^{n+1}N} = 0 = \tor iRM{\fm^{n+1}N}
\]
then either $\fm^{n+1}N=0$ or $\pd_RM\leq i-1$.
\end{theorem}

These theorems are akin to the one below, contained in \cite[Lemma, p. 316]{LV}.

\begin{theorem}
\label{theorem:lv}
Let $(R,\fm,k)$ be a local ring, and let $M$ and $N$ be  finite $R$-modules. 
If there exist non-negative integers $i$ and $n$ such that
\[
\tor iRM{\fm^{n+1}N} = 0 = \tor {i+1}RM{\fm^{n+1}N}
\]  
then either $\fm^{n+1}N=0$ or $\pd_RM\leq i-1$. \qed
\end{theorem}

Evidently, Theorems \eqref{Intheorem}, \eqref{Intheorem2}, and \eqref{theorem:lv} are
closely related: consider the exact sequence of $R$-modules $0\to \fm^{n+1}N\to N\to
N/\fm^{n+1}N\to 0$. However, we have been unable to deduce one from the other, nor find a
useful common generalization.  We should like to note that statements analogous to the
ones above with $\tor *RM-$ replaced by $\ext *RM-$ also hold.  Now we turn to the

\begin{proof}[Proof of Theorems  \eqref{Intheorem}  and \eqref{Intheorem2}]
  We may assume that $M$ and $N$ are non-zero. Thus, $\tor 0RMN = M\otimes_RN\ne 0$, so
  $i\geq 1$, and
\[
\tor iRM- \cong \tor 1R{\syz {i-1}RM}-\quad\text{and}\quad
\pd_RM = \pd_R{\syz iRM} + i
\]
Hence, replacing $M$ with $\syz{i-1}RM$ one may assume that $i=1$. Set $L=\syz 1RM$ and
$F=R^{\nu(M)}$, and consider the  exact sequence of $R$-modules
\begin{equation*}
0\to L\to F\to M\to 0 \tag{$\dagger$}
\end{equation*}
\emph{Claim}:
There is a commutative diagram of homomorphisms of $R$-modules:
\[
\xymatrixrowsep{1.5pc} 
\xymatrixcolsep{2.5pc} 
\xymatrix{ 
               & 0\ar@{->}[d] & 0\ar@{->}[d] & 0\ar@{->}[d]  \\
0\ar@{->}[r]   & \Ker(\alpha)\ar@{->}[r] \ar@{->}[d]& L\otimes_R N \ar@{->}[r]^-{\beta} \ar@{->}[d]
               & L\otimes_R (N/\fm^{n+1}N) \ar@{->}[r] \ar@{->}[d] & 0 \\
0 \ar@{->}[r]  & F\otimes_R \fm^{n+1}N \ar@{->}[r] \ar@{->}[d]^-{\alpha}
               & F\otimes_R N \ar@{->}[r] \ar@{->}[d]
               & F\otimes_R (N/\fm^{n+1}N) \ar@{->}[r] \ar@{->}[d] & 0 \\
0 \ar@{->}[r]  & M\otimes_R \fm^{n+1}N \ar@{->}[r] 
               & M\otimes_R N \ar@{->}[r] 
               & M\otimes_R (N/\fm^{n+1}N) \ar@{->}[r] & 0 
}
\]
where the sequences in the rows and columns are exact, and, under the hypotheses of
\eqref{Intheorem2}, also $\Ker(\alpha)=L\otimes_R\fm^{n+1}N$.

Indeed, tensoring ($\dagger$) with the exact sequence 
\begin{equation*}
 0\to \fm^{n+1}N \to N \to N/\fm^{n+1} N\to 0 \tag{$\ddagger$}
\end{equation*}
yields the desired commutative diagram. Now, under either sets of hypotheses, $\tor
1RM{N/\fm^{n+1}N}=0$, so the last row and the rightmost column are also exact.  The second
row is exact because the $R$-module $F$ is free, while the second column is exact because
$\tor 1RMN=0$; this is part of the hypotheses in \eqref{Intheorem}, and follows from that
of \eqref{Intheorem2} and the exact sequence ($\ddagger$) above.  The exactness of the
last two rows and the surjectivity of $\alpha$ imply, by the snake lemma, that the first
row is exact.  Finally, in case \eqref{Intheorem2}, the leftmost column is exact and
$\Ker(\alpha)=L\otimes_R\fm^{n+1}N$.

To finish the proof, we consider the two cases separately. In either case we may suppose
that $\pd_RM\geq 1$, that is to say, $M$ is not free. Thus, $L\ne 0$.

\smallskip

\emph{Case \eqref{Intheorem}}. Viewing $F\otimes_R\fm^{n+1}N$ and $L\otimes_RN$ as submodules of
$F\otimes_RN$, one has $\Ker(\alpha)=\Ker(\beta)$.  This explains the second equality
below
\[
\left(L\otimes_R N\right)\cap \fm^{n+1} \left(F\otimes_R N\right)
  = \left(L\otimes_R N\right)\cap \left(F\otimes_R \fm^{n+1} N\right)
  = \fm^{n+1}\left( L\otimes_R N\right)
\]
while the first is due to the $R$-linearity of the tensor product.  On the other hand,
$L\subseteq \fm F$, so $L\otimes_R N \subseteq \fm\left( F\otimes_R N\right)$, and hence
$\fm^n \left(L\otimes_RN\right) \subseteq \fm^{n+1}\left( F\otimes_R N\right)$.  Combining
this inclusion with the equality above yields
\[
\fm^n \left(L\otimes_R N\right) \subseteq \left(L\otimes_R N\right)\cap
\fm^{n+1}\left(F\otimes_R N\right) = \fm^{n+1}\left( L\otimes_R N\right) .
\]
Thus $\fm^n\left(L\otimes_RN\right) = 0$, by Nakayama's lemma. This is the first part of
the desired result.  As to the second, since $L$ is non-zero $L\otimes_RN$ is non-zero as
well.  Since $(L\otimes_RN)\subseteq (F\otimes_R N)$, by the diagram above, and
$\fm^n\left(L\otimes_RN\right)=0$, we deduce that $\depth (F\otimes_RN)=0$. However,
$F\otimes_R N\cong N^{\nu(M)}$, so $\depth N=0$, as claimed.

\smallskip

\emph{Case \eqref{Intheorem2}}.  Recall that $\Ker(\alpha)= L\otimes_R\fm^{n+1}N$.
Viewing the modules in the top left-hand square of the diagram above as submodules of
$F\otimes_RN$, one has that \( L\otimes_R\fm^{n+1}N = \fm^{n+1}\left(L\otimes_RN\right) =
0 \) where the second equality holds because $\fm^n\left(L\otimes_RN\right) = 0$, by case
\eqref{Intheorem}. However, $L\ne 0$, so $\fm^{n+1}N=0$, as desired.
\end{proof}

Here is one consequence of Theorem \eqref{Intheorem2}.

\begin{corollary}
\label{testmodule}
Let $(R,\fm,k)$ be a local ring, $N$ a finite $R$-module with
$\fm^{n+1}N\ne 0$ and $\pd_R(N/\fm^{n+1}N)$ finite for an integer $n\ge0$.
For each finite $R$-module $M$ with $\pd_RM=\infty$, one has $\tor iRMN\ne 0$ for each
integer $i \ge\depth R+1$.
\end{corollary}
\begin{proof}
  Fix an integer $i\ge\depth R +1$.  The Auslander-Buchsbaum equality implies
  $\pd_R(N/\fm^{n+1}N)=\depth R$, so $\tor jRM{N/\fm^{n+1}N}=0$ for $j=i,i+1$. Thus, for
  each $R$-module $M$ the exact sequence $0\to \fm^{n+1}N\to N\to N/\fm^{n+1}N\to 0$
  yields $\tor iRM{\fm^{n+1}N}\cong \tor iRMN$. Therefore, the desired result follows from
  (the contrapositive of) Theorem \eqref{Intheorem2}.
\end{proof}

Next we introduce some invariants used to describe other results in this section.

\begin{chunk}
\label{avind:definition}
  Let $N$ be a finite $R$-module. For each non-negative integer $n$, let $\pi^n$ denote
  the canonical surjection $N\to N/\fm^{n+1}N$. We define the \emph{Avramov index} of
  $N$ to be 
\[
\avind_R(N) = \inf\{n\geq 0 \mid \tor{}Rk{\pi^n}\quad\text{is injective.}\}
\]
This definition is motivated by a result of Avramov \cite[(A.4)]{Av1}, which implies that
$\avind_R(N)$ is finite.  \c{S}ega \cite[(5.1)]{Sega} has introduced the Avramov index of
a local \emph{ring} $R$ to the least integer $n$ such that the map
$\tor{}Rkk\to\tor{}{R/\fm^{n+1}}kk$, induced by the canonical homomorphism of rings $R\to
R/\fm^{n+1}$, is injective. It is not hard to prove that $\avind_R(\fm)$ is less than or
equal to the Avramov index of the ring $R$; we do not know if equality holds.

The \emph{Levin index} of $N$, see \cite[(3.1)]{Sega}, is the least integer $n$ for which
the map $\tor{}Rk{\fm^s N}\to\tor{}{R}k{\fm^{s-1}N}$, induced by $\fm^s N\subseteq
\fm^{s-1}N$, is zero for $s\geq n$.
\end{chunk}

The invariant $\preg_R(N)$ appearing in the lemma below is defined in \eqref{polyreg}.

\begin{lemma}
\label{avind}
Let $R$ be a local ring. For each finite $R$-module $N$, one has
\[
\avind_R(N)\leq \lind_R(N) - 1 \leq \preg_R(N)
\]
\end{lemma}

\begin{proof}
Applying $k\otimes_R-$ to the exact sequence $0\to \fm^{n+1} N\to N\to N/\fm^{n+1} N\to 0$
yields the long exact sequence of $R$-modules
\[
\cdots \to \tor {i}Rk{\fm^{n+1} N} \to \tor {i}RkN \xra{\tor iRk{\pi^n}}
    \tor {i}Rk{N/\fm^{n+1} N} \to \cdots
\]
Thus, $\tor {}Rk{\pi^n}$ is injective if and only if the map $\tor {}Rk{\fm^{n+1}
  N}\to\tor{}RkN$ is zero. Since this last map factors through $\tor{}Rk{\fm^{n}N}$,
one obtains that if $\tor{}Rk{\fm^{n+1} N}\to\tor{}Rk{\fm^{n} N}$ is zero, then
$\tor{}Rk{\pi^n}$ is injective as well. Thus, $\avind_R(N)\leq \lind_R(N)-1$. The remaining
inequality is contained in  \cite[(3.3)]{Sega}.
\end{proof}

The next result is easily deduced from available literature, but there is no convenient
reference; in any case, it is worth stating it in terms of the Avramov index.

\begin{theorem}
\label{nth}
Let $(R,\fm,k)$ be a local ring. If there is a finite $R$-module $N$ with $\dim_RN\ge
1$ and $\pd_{R}(N/\fm^{n+1}N)$ finite for an  $n\geq \avind_R(N)$, then $R$ is regular
\end{theorem}

\begin{proof}
The condition $n\geq \avind_R(N)$ ensures the injectivity of the map
\[
\Tor^R(k,\pi^n)\colon \Tor^R(k,N)\lra \Tor^R(k,N/\fm^{n+1}N)
\]
Therefore, since $\pd_R(N/\fm^{n+1}N)$ is finite, $\Tor^R_i(k,N)=0$ for $i\gg 0$, so that
$\pd_R(N)$ is also finite. If $\depth_RN\geq1$, then at this point Corollary
\eqref{testmodule} applied with $M=k$ would entail $\pd_Rk$ finite, and hence that $R$ is
regular.  In any case, from the exact sequence
\[
0\to \fm^{n+1}N\to N\to N/\fm^{n+1}N\to 0
\]
one deduces that $\pd_R(\fm^{n+1}N)$ is finite. Since $\dim_RN\ge 1$, Nakayama's
lemma implies  $\fm^{n+1}N\ne 0$. It now follows from \cite[(1.1)]{LV} that $R$ is regular.
\end{proof}

The preceding theorem has an analogue for injective dimensions, where, instead of the
Avramov index, the invariant of interest is the least integer $n$ such that the canonical
homomorphism $\Ext_R(k,N)\to\Ext_R(k,N/\fm^{n+1}N)$ is injective. We take a different
route in the ensuing result; the invariant $\rho_R(N)$ that appears here was introduced in
\eqref{rho vs polyreg}. As noted there, it is bounded above by $\preg_R(N)$.
Recall that $R$ is said to be a \emph{hypersurface} if $\edim R - \depth R\leq 1$.

\begin{theorem}
\label{nth2}
Let $(R,\fm,k)$ be a local ring and let $N$ be a finite $R$-module with $\depth_RN\geq 1$
and $\id_{R}(N/\fm^{n+1}N)$ finite for some non-negative integer $n$.

If $n \geq \rho_R(N)$, then $R$ is a hypersurface; if $n\geq \avind_R(N)$ as well, then
$R$ is regular.
\end{theorem}

\begin{proof}
  Let $L$ be an $R$-module; we set $s_i(L) =\length_R\, \Hom_R(L, N/\fm^{i+1}N)$ and
  write $\nu_R(L)$ for the minimal number of generators $L$.  First we establish the
claims below.

\smallskip

\emph{Claim}.  One has $s_n(L) \geq s_{n-1}(L)$.

Indeed, \eqref{infresfield} allows one to enlarge the residue field and
assume that the field $k$ is infinite. Thus, there is an $x$ in $\fm$
superficial on $N$, so for $n \geq \rho(N)$ multiplication by $x$
yields an exact sequence of $R$-modules \( 0 \lra N/\fm^{n}N \to N/\fm^{n+1}N, \) and this
gives the exact sequence below, which settles the claim:
\[
0 \lra \Hom_R(L ,N/\fm^{n}N) \lra \Hom_R(L, N/\fm^{n+1}N ).
 \]

\smallskip

\emph{Claim}.
  If $\Ext_{R}^{1}(L, N/\fm^{n+1}N) = 0$, then $\nu_R(L) \geq \nu_R \left(\syz 1RL\right)$.
 
Indeed, the exact sequence of $R$-modules
\[
0 \lra {\fm^nN}/{\fm^{n+1}N} \lra N/\fm^{n+1}N \lra N/\fm^nN \lra 0
\]
induces an exact sequence 
\begin{align*}
  0 &\lra \Hom_R(L, {\fm^nN}/{\fm^{n+1}N} ) 
            \lra  \Hom_R(L, N/\fm^{n+1}N) \lra  \Hom_R(L, N/\fm^{n}N)  \\
  &\lra \Ext^{1}_{R}(L, {\fm^nN}/{\fm^{n+1}N} ) \lra \Ext^{1}_{R}(L,N/\fm^{n+1}N) = 0
\end{align*}
Note that $\Ext^{1}_{R}(L, {\fm^nN}/{\fm^{n+1}N})\cong \Hom _R(\syz 1RL,
{\fm^nN}/{\fm^{n+1}N})$, so computing lengths one obtains from the exact sequences above, an
equality
\[
h\cdot \nu_R (L) - h\cdot \nu_R(\syz 1RL) = s_n(L) - s_{n-1}(L)
\]
where $h=\rank_k(\fm^n N/\fm^{n+1}N)$. Since $\dim_R N \geq 1$ one has $h\geq 1$, so the
equality above and the preceding claim imply that $\nu_R(L)\geq \nu_R(\syz 1RL)$, as desired.

Consider the case when $d=\id_R(N/\fm^{n+1}N)$ is finite. For $i\geq d+1$, one has
\[
\Ext_R^1(\syz iRk,N/\fm^{n+1}N) = \Ext_R^{i+1}(k,N/\fm^{n+1}N)=0
\]
Given this, the preceding claim yields $\nu_R(\syz iRk)\geq \nu_R(\syz {i+1}Rk)$.  Thus,
the Betti numbers of $k$ are bounded, so $R$ is a hypersurface; see Tate \cite[Theorem
8]{Tate}. In particular, $R$ is a hypersurface, and hence Gorenstein, so
$\pd_R(N/\fm^{n+1}N)$ is finite. As $\depth_RN\geq1$, when $n\geq \avind_R(N)$ Theorem
\eqref{nth} implies $R$ is regular.
\end{proof}

This last result leads us to modules of finite length and finite projective dimension over
hypersurfaces. In this context one has the following well known result; see, for example,
Ding \cite[(1.5) and (3.3)]{Di}.

\begin{proposition}
\label{hyp}
Let $(R,\fm,k)$ be a hypersurface and $M$ a non-zero finite $R$-module.
If $\fm^{e(R)-1}M=0$, then 
\(
\pd_RM =\infty = \id_R M
\)
\end{proposition}
\begin{proof}
  It is elementary to reduce to the case where $R$ and $M$ are $\fm$-adically complete.
  Now, Cohen's structure theorem provides a surjective homomorphism $(Q,\fq,k)\to R$ of
  local rings with $Q$ regular and $\edim Q=\edim R$. Since $R$ is a hypersurface,
  $R=Q/(r)$; it is not hard to verify that $r\in\fq^{e(R)}$. Note that $\fq^{e(R)-1}M=0$,
  so $r\in \fq \Ann_Q(M)$, where $\Ann_Q(M)$ denotes the annihilator of $M$ viewed as a
  module over $Q$. A result of Shamash \cite[Corollary 1]{Sh} now implies that $\pd_RM$ is
  infinite.  Since $R$ is a Gorenstein, $\id_RM$ is also infinite.
\end{proof}

>From Theorem \eqref{nth2} and Proposition \eqref{hyp} one obtains the following corollary.
A hypersurface $R$ that is not regular is a \emph{singular hypersurface}; thus, $R$ is a
singular hypersurface when $\edim R - \depth R =1$.

\begin{corollary}
  Let $R$ be a singular hypersurface. Let $N$ be a finite $R$-module with $\depth_RN\geq1$
  and set $d=\sup\{\rho_R(N),\avind_R(N)\}$. For each integer $n\ge0$ with
  $n\not\in [e(R)-1,d-1]$, one has
\(
\pd_R(N/\fm^{n+1}N)=\infty=\id_R(N/\fm^{n+1}N)
\)\qed
\end{corollary}




We end with a question: does the conclusion of the corollary hold for each $n\geq0$?

\section*{Acknowledgments}
It is a pleasure to thank Lucho Avramov, Shiro Goto, Janet Striuli, and Emanoil
Theodorescu for useful conversations regarding this article.

\end{document}